\documentclass[a4paper]{amsart}
\usepackage{latexsym}
\usepackage{amsmath}
\usepackage{amsfonts}
\usepackage{mathrsfs}
\usepackage{amscd}
\usepackage{amssymb}
\usepackage{amsfonts}
\usepackage{float}
\usepackage{esint}
\usepackage{caption}
\usepackage[utf8]{inputenc}
\usepackage[T1]{fontenc}   
\usepackage{graphicx}
\usepackage{color}
\allowdisplaybreaks 
\definecolor{Blue}{rgb}{0.,0.,1.}
\definecolor{Red}{rgb}{1.,0.,0.}
\definecolor{Green}{rgb}{0.,1.,0.}
\let\origmaketitle\maketitle
\def\maketitle{
  \begingroup
  \def\uppercasenonmath##1{} 
  \let\MakeUppercase\relax 
	\origmaketitle
  \endgroup
	}

\newcounter{smallarabics}
\newenvironment{arabicenumerate}
{\begin{list}{{\normalfont\textrm{(\arabic{smallarabics})}}}
  {\usecounter{smallarabics}\setlength{\itemindent}{0cm}
   \setlength{\leftmargin}{5ex}\setlength{\labelwidth}{4ex}
   \setlength{\topsep}{0.75\parsep}\setlength{\partopsep}{0ex}
   \setlength{\itemsep}{0ex}}}
{\end{list}}

\newcounter{smallroman}

\newcommand{\ben}{\begin{arabicenumerate}}  
\newcommand{\een}{\end{arabicenumerate}}

\def\init{\setcounter{equation}{0}}

\newtheorem{theoreme}{Theorem }[section]
\newtheorem{proposition}[theoreme]{Proposition}

\newtheorem{lemma}[theoreme]{Lemma}
\newtheorem{definition}[theoreme]{Definition}

\newtheorem{remark}[theoreme]{Remark}
\newtheorem{example}[theoreme]{Example}
\newcommand{\beq}{\begin{equation}}
\newcommand{\eeq}{\end{equation}}

\newcommand{\bex}{\begin{example}}
\newcommand{\eex}{\end{example}}
\def\bel{\begin{lemma}}
\def\eel{\end{lemma}}
\def\bet{\begin{theoreme}}
\def\eet{\end{theoreme}}
\def\bed{\begin{definition}}
\def\eed{\end{definition}}
\def\ber{\begin{remark}}
\def\eer{\end{remark}}

\def\rr{{\mathbb R}}

\def\cc{{\mathbb C}}
\def\nn{{\mathbb N}}

\def\bar{\overline}

\def\cinf{C^\infty}
\def\proof{
\noindent{\bf Proof.}\ \ }

\def\cY{{\pazocal Y}}

\def\cD{{\mathcal D}}
\def\cU{{\mathcal U}}

\def\i{{\rm i}}
\def\qed{$\Box$\medskip}
\def \p{ \partial}
\def\12{\frac{1}{2}}
\def\14{\frac{1}{4}}

\def\bbbone{{\mathchoice {\rm 1\mskip-4mu l} {\rm 1\mskip-4mu l}
{\rm 1\mskip-4.5mu l} {\rm 1\mskip-5mu l}}}
\def\one{\bbbone}
\def\cH{{\mathcal H}}

\def\coinf{C_0^\infty}

\def\cX{{\pazocal X}}

\def \p{ \partial}
\def\12{\frac{1}{2}}
\def\e{{\rm e}}

\def\d{{\rm d}}

\newcommand{\mat}[4]{\left(\begin{array}{cc}#1 &#2  \\ #3 &#4 \end{array}\right)}

\newcommand{\col}[2]{\left(\begin{array}{c}#1 \\#2\end{array} \right)}

\newcommand{\traa}[1]{\mskip-6mu\upharpoonright_{#1}}

\def\cE{{\mathcal E}}

\def\WF{{\rm WF}}
\makeatletter
\newcommand*{\defeq}{\mathrel{\rlap{%
                     \raisebox{0.3ex}{$\m@th\cdot$}}%
                     \raisebox{-0.3ex}{$\m@th\cdot$}}%
                     =}
\makeatother
\makeatletter
\newcommand*{\eqdef}{=\mathrel{\rlap{%
                     \raisebox{0.3ex}{$\m@th\cdot$}}%
                     \raisebox{-0.3ex}{$\m@th\cdot$}}%
                     }
\makeatother

\def\rx{{\rm x}}

\DeclareMathOperator{\Ker}{Ker}
\DeclareMathOperator{\Dom}{Dom}
\DeclareMathOperator{\supp}{supp}

\def\dual{\!\cdot \!}

\def\zero{{\mskip-4mu{\rm\textit{o}}}}

\def\cN{{\mathcal N}}

\def\outin{{\rm out/in}}

\def\adg{{\rm ad}}

\def\pe{\p}
\def\std{{\rm std}}
\def\free{{\rm free}}

\def\altg{{\rm\textit{g}}}
\def\alth{{\rm\textit{h}}}
\def\altV{{\rm\textit{V}}}

\def\altm{{\rm\textit{m}}}

\newcommand{\bea}{\begin{aligned}}
\newcommand{\beal}{\begin{array}{l}}
\newcommand{\eeal}{\end{array}}
\newcommand{\eea}{\end{aligned}}
\def\sobo{{m}}

\def\scc{{\rm sd}}
\def\spexi{{k}}
\usepackage{calrsfs}
\DeclareMathAlphabet{\pazocal}{OMS}{zplm}{m}{n}
\DeclareMathAlphabet{\mathsfsl}{OMS}{cmss}{m}{n}

\DeclareSymbolFont{altletters}  {OML}{zplm}{m}{it}
\DeclareMathSymbol{\altdelta}{\mathalpha}{altletters}{"0E}
\DeclareMathSymbol{\alteta}{\mathalpha}{altletters}{"11}

\def\altV{{\rm\textit{V}}}
\def\F{{\rm F}}
\def\diag{{\rm diag}}
\begin{document}
\pagestyle{plain}
\title[]{\Large The massive Feynman propagator \\ on asymptotically  Minkowski spacetimes II}
\address{Universit\'e Paris-Sud XI, D\'epartement de Math\'ematiques, 91405 Orsay Cedex, France}
\email{christian.gerard@math.u-psud.fr}
\author{\large Christian G\'erard  \& Micha{\l} Wrochna}
\thanks{\emph{Acknowledgments.} We would like to warmly thank Christian B\"{a}r, Jan Derezi\'nski, Alexander Strohmaier and Andr\'as Vasy for useful conversations, and in particular for pointing out to us the link between invertibility and general features of the Klein-Gordon operator. Support from the grant ANR-16-CE40-0012-01 is gratefully acknowledged.}
\address{Universit\'e Grenoble Alpes, Institut Fourier, UMR 5582 CNRS, CS 40700, 38058 Grenoble \textsc{Cedex} 09, France}
\email{michal.wrochna@univ-grenoble-alpes.fr}
\begin{abstract}
We consider the massive Klein-Gordon equation on short-range asymptotically Minkowski spacetimes. Extending our results in \cite{GW1}, we show that the Klein-Gordon operator with Feynman type boundary conditions at infinite times is invertible and that its inverse, called the {\em Feynman inverse},  satisfies the microlocal conditions of Feynman parametrices 
 in the sense of Duistermaat and H\"ormander. This supplements the recent work of Vasy \cite{Va} with more explicit techniques.
\end{abstract}

\maketitle

\section{Introduction}\label{sec0}\init
The present paper is a continuation of \cite{GW1}, which was devoted to the existence of the  {Feynman propagator} for Klein-Gordon fields on asymptotically Minkowski spacetimes. Let us first briefly recall the motivation of \cite{GW1}. 

On Minkowski spacetime $(\rr^{1+d}, \alteta)$, the Klein-Gordon operator $P_{\free}= \p_{t}^{2}- \Delta_{\rx}+\altm^{2}$ has four distinguished inverses, the {\em retarded/advanced} inverses $G_{\rm ret/adv}$ which are Fourier multipliers by $((\tau\pm\i 0)^{2}- (k^{2}+\altm^{2}))^{-1}$, and the {\em Feynman/anti-Feynman} inverses $G_{{\rm F}/\overline{\rm F}}$ which are Fourier multipliers by $((\tau)^{2}- (k^{2}+\altm^{2})\pm \i 0)^{-1}$.

The retarded/advanced inverses exist on any globally hyperbolic spacetime $(M, g)$ and are characterized as the unique solutions of 
\[
\begin{array}{l}
PG_{\rm ret/adv}= G_{\rm ret/adv}P= \one, \\[2mm]
 \supp G_{\rm ret/adv}v\subset J^{\pm}\left(\supp v\right), \ v\in \coinf(M),
\end{array}
\]
where $P= - \Box_{g}+ V$, $V\in\cinf(M; \rr)$ and $J^{\pm}(K)$ is the future/past causal shadow of a set $K\subset M$.  Their difference $G= G_{\rm ret}- G_{\rm adv}$, called the {\em causal} or {\em Pauli-Jordan propagator} is used in  the algebraic quantization of free Klein-Gordon fields on $(M, g)$.

The Feynman/anti-Feynman inverses on Minkowski spacetime play a  fundamental role in the perturbative renormalization of interacting Klein-Gordon fields.

It is not a priori clear  how the  Feynman/anti-Feynman inverses generalize to an arbitrary globally hyperbolic spacetime $(M, g)$, but some of the properties they need to satisfy  are understood since the work of Duistermaat and H\"ormander \cite{DH}, who proved the existence and uniqueness modulo smooth terms of \emph{Feynman parametrices}.

Namely, if  $\Phi_t$ is the Hamilton flow of $p(x,\xi)=\xi\cdot \altg^{-1}(x)\xi$ restricted to the \emph{characteristic set} $\cN=p^{-1}(\{0\})$ (understood as a subset of $T^*M\setminus\zero$, where $\,\zero$ is the zero section of the cotangent bundle), one says that 
 $\tilde{G}_\F$ is a {\em Feynman parametrix} if the operators $\one- \tilde{G}_\F P$ and $\one-P \tilde{G}_\F$ have smooth Schwartz kernels and
\beq\label{eq:fewf}
\WF'(\tilde{G}_{\F})= (\diag_{T^*M})\cup\textstyle\bigcup_{s\leq 0}(\Phi_s(\diag_{T^*M})\cap \pi^{-1}\cN).
\eeq
Here  $\Phi_s$ is the bicharacteristic flow acting on the left component of $\diag_{T^*M}$ (i.e., the diagonal in $(T^*M\times T^*M)\setminus\zero$), and $\pi:\cN\times\cN\to\cN$ is the projection to the left component.

This leaves open the question of the existence of a {\em canonical Feynman inverse} $G_{\F}$ satisfying \eqref{eq:fewf} on a globally hyperbolic spacetime $(M, g)$, see \cite{DS2} for a discussion. 

If $(M, g)$ is {\em stationary}, i.e.~admits a global, complete, time-like Killing vector field $K$ and if  the potential  $V$ is invariant under $K$, it is well known that the {\em vacuum state} (with respect to the group of isometries generated by $K$) generates a canonical Feynman inverse.

In \cite{GW1} we considered this problem for {\em asymptotically Minkowski} spacetimes. Inspired by  works by Gell-Redman, Haber and Vasy \cite{GHV,positive} and B\"{a}r and Strohmaier \cite{BS} on closely related non-elliptic Fredholm problems, we introduced Hilbert spaces $\cY^{m}$, $\cX^{m}_{\F}$, where $m\in \rr$ is an order of Sobolev regularity and the subscript $\F$ refers to what one can call {\em Feynman boundary conditions}, which are imposed at $t= \pm \infty$.  We refer the reader to \ref{sec0.1.1} and \ref{sec0.2.1} for precise definitions.

We proved in \cite[Thm.~1.2]{GW1} in the massive case that $P:\cX^{m}_{\F}\to \cY^{m}$ is a {\em Fredholm operator} and we related  its index  to the index of  some `Feynman' wave operator $W_{\F}^{\dag}$.  We also showed that $P$ has a pseudo inverse $\tilde{G}_{\F}: \cY^{m}\to \cX^{m}_{\F}$ (i.e., an inverse modulo compact errors) which is at the same time a Feynman parametrix. Though obtained with different methods, these results are analogous to the outcome of \cite{GHV,VW} in the massless case, and are consistent with the general program outlined in Vasy's work \cite{positive}.

Our method of proof consisted in writing the Klein-Gordon equation as a first order system, which can be diagonalized globally in time modulo smoothing and decaying in time errors. After these reductions the problem can be handled by rather explicit methods.

Recently, Vasy \cite{Va} considered the same problem by working directly on the scalar operator $P$ using microlocal methods, motivated by the issue of essential self-adjointness of $P$ (see \cite{DS1}). Solving in this setting a conjecture by Derezi\'nski and Siemssen \cite{DS2,D}, he constructed the Feynman inverse $G_{\F}$ between microlocal Sobolev spaces as the boundary value $(P-\i 0)^{-1}$ of the resolvent of $P$. 

In the present paper we  will  recover Vasy's result in our framework by showing  that $P: \cX_{\F}^{m}\to \cY^{m}$ is indeed {\em invertible} rather than merely Fredholm, and that its inverse $G_{\F}$ is a {\em Feynman inverse}, i.e.~the wavefront set of its kernel equals the r.h.s.~in \eqref{eq:fewf},  see Thm.~\ref{thm1} below.

The fact that $P$ is of index zero could actually be concluded directly from \cite{GW1} (see the proof of Prop.~\ref{prop1.1} in the present paper). Its injectivity from $\cX^{m}_{\F}$ to $\cY^{m}$ is less evident, the proof turns out however to be very easy, by adapting arguments of \cite{Va} to our framework. The advantage of having such a result in our framework is that this provides a rather explicit parametrix for $G_{\rm F}$.

\subsection{Klein-Gordon operators on asymptotically Minkowski spacetimes}\label{sec0.1}
In this subsection we recall the framework considered in \cite{GW1}.
\subsubsection{Asymptotically Minkowski spacetimes}\label{sec0.1.1}
We consider    $M=\rr^{1+d}$ equipped with  a Lorentzian metric $\altg$ such that
\[
 \begin{array}{rl}
({\rm aM} i)&\altg_{\mu\nu}(y)- \alteta_{\mu\nu} \in S^{-\delta}_{\std}(\rr^{1+d}), \ \delta>1,\\[2mm]
({\rm aM} ii)&(\rr^{1+d}, \altg) \hbox{ is globally hyperbolic},\\[2mm]
({\rm aM} iii)&(\rr^{1+d}, \altg) \hbox{ has a  time function }\tilde{t}\hbox{ with  }\tilde{t}- t\in S^{1-\epsilon}_{\std}(\rr^{1+d})\hbox{ for }\epsilon>0.
\end{array}
\]
where $\alteta_{\mu\nu}$ is the Minkowski metric and $S_{\std}^{\delta}(\rr^{1+d})$ stands for the class of smooth functions $f$ such that for  $\langle x \rangle  = (1+|x|)^{\12}$,
\[
\p^{\alpha}_{x}f\in O(\langle x\rangle^{\delta- |\alpha|}), \ \alpha\in \nn^{1+d}.
\] 
We recall $\tilde{t}$ is a {\em time function} if $\nabla \tilde{t}$ is a time-like vector field.  It is a {\em Cauchy time function} if in addition its level sets  are Cauchy surfaces for $(M, \altg)$. 

It is shown in \cite{GW1} that if $({\rm aM} i)$ holds, then $({\rm aM} ii)$ is equivalent to the familiar {\em non trapping condition} for null geodesics of $\altg$, and if $({\rm aM} i), ii) , iii)$ hold there exists a  Cauchy time function $\tilde{t}$ such that  $\tilde{t}- t\in \coinf(M)$.

Replacing $t$ by $t-c$, $\tilde{t}$ by $\tilde{t}-c$ for $c\gg 1$ we can also assume that $\Sigma\defeq \{t=0\}= \{\tilde{t}=0\}$ is a Cauchy surface for $(M, \altg)$, which can be canonically identified with $\rr^{d}$. In the sequel we will fix such a time function $\tilde{t}$.

\subsubsection{Klein-Gordon operator}\label{sec0.1.2}
We  fix a real function $V\in \cinf(M; \rr^{d})$ such that
\[
({\rm aM} iv) \quad V(x)- \altm^{2} \in S^{-\delta}_{\std}(\rr^{1+d}),\hbox{ for some }\altm>0,\ \delta>1,
\]
and consider the Klein-Gordon operator
\[
P= - \Box_{\altg}+V.
\]
\subsection{The Feynman inverse of $P$}\label{sec0.2}
We now introduce the Hilbert spaces  $\cX^{m}_{\F}$, $\cY^{m}$ between which $P$ will be  invertible. The spaces $\cY^{m}$ are standard spaces of right hand sides for the Klein-Gordon equations, their essential property being that they are $L^{1}$ in $t$, with values in some Sobolev spaces of order $m$ (not to be confused with the Klein-Gordon mass parameter $\altm$). The spaces $\cX^{m}_{\F}$ incorporate the {\em Feynman boundary conditions}, which are imposed at $t= \pm \infty$.

\subsubsection{Hilbert spaces}\label{sec0.2.1}
Using the Cauchy time function $\tilde{t}$ we can  identify $M$ with $\rr\times \Sigma$, using the flow $\phi_{t}$ of the vector field $v= \frac{\altg^{-1}d\tilde{t}}{d\tilde{t}\cdot \altg^{-1} d\tilde{t}}$ and obtain  the diffeomorphism:
\beq\label{e0.10}
\chi: \rr\times \Sigma\ni (t, \rx)\to \phi_{t}(\rx)\in M,
\eeq
such that
\[
\chi^{*}\altg= - c^{2}(t, \rx)dt^{2}+ h(t, \rx)d\rx^{2}.
\]
For $m\in \rr$ we denote by $H^{m}(\rr^{d})$ the usual Sobolev spaces on $\rr^{d}$. We set for $\12<\gamma<\12+ \delta$:
\[
\cY^{m}\defeq \{u\in \cD'(M): \chi^{*}u\in \langle t\rangle^{-\gamma}L^{2}(\rr; H^{m}(\rr^{d}))\}
\]
with norm $\|v\|_{\cY^{m}}= \| \langle t\rangle^{\gamma}\chi^{*}u\|_{L^{2}(\rr; H^{m}(\rr^{d}))}$.  The exponent $\gamma$ is chosen such that $ \langle t\rangle^{-\gamma}L^{2}(\rr)\subset L^{1}(\rr)$. Similarly we set
\[
\cX^{m}\defeq \{u\in \cD'(M): \chi^{*}u\in C^{0}(\rr; H^{m+1}(\rr^{d}))\cap C^{1}(\rr; H^{m}(\rr^{d})), \ Pu\in \cY^{m}\}.
\]
We equip  $\cX^{m}$ with the norm
\[
\| u\|_{\cX^{m}}= \| \varrho_{0} u\|_{\cE^{m}}+ \| Pu\|_{\cY^{m}},
\]
where $\varrho_{s}u= \col{u\traa{\Sigma_{s}}}{\i^{-1}\p_{n}u\traa{\Sigma_{s}}}$ is the Cauchy data map on $\Sigma_{s}\defeq \tilde{t}^{-1}(\{s\})$  and $\cE^{m}\defeq H^{m+1}(\rr^{d})\oplus H^{m}(\rr^{d})$ is the  energy space of order $m$. From the well-posedness of the inhomogeneous Cauchy problem for $P$ one easily deduces that $\cX^{m}$ is a Hilbert space. 

\subsubsection{Feynman boundary conditions}\label{sec0.2.2}
Let us set 
\[
c^{\pm}_{\free}=\12\begin{pmatrix}1 & \pm \sqrt{-\Delta_\rx+\altm^2} \\ \pm \sqrt{-\Delta_\rx+\altm^2} & 1\end{pmatrix},
\] 
which are the spectral projections on $\rr^{\pm}$ for the generator 
\[
H_{\free}= \mat{0}{1}{-\Delta_\rx+\altm^2}{0}
\]
of the Cauchy evolution for the free Klein-Gordon operator $P_{\free}= \p_{t}^{2}- \Delta_{\rx}+ \altm^{2}$, to which $P$  is asymptotic when $t\to \pm \infty$.

We then set
\[
\cX^{m}_{\rm F}\defeq \{u\in \cX^{m}: \lim_{t\to \mp\infty}c^{\pm}_{\free}\varrho_{t}u=0\hbox{ in }\cE^{m}\}.
\]
It is easy to see that $\cX^{m}_{\F}$ is a closed subspace of $\cX^{m}$.

In this paper we will prove the following theorem:
\begin{theoreme}\label{thm1}
 Assume $({\rm aM})$. Then $P: \cX^{m}_{\rm F}\to \cY^{m}$ is boundedly invertible for all $m\in\rr$. Its inverse $G_{\rm F}$ is called the {\em Feynman inverse} of $P$. It satisfies:
 \[
  \WF(G_{{\rm F}})'= ({\rm diag}_{T^*M})\cup\textstyle\bigcup_{t\leq 0}(\Phi_t({\rm diag}_{T^*M})\cap \pi^{-1}\cN).
  \]
\end{theoreme}
\section{Proof of Thm.~\ref{thm1}}\label{sec1}\init
We now give the proof of Thm.~\ref{thm1}. The first step consists in replacing $P$ by $\chi^{*}P$, where $\chi$ is the diffeomorphism in \eqref{e0.10}. After an additional conformal transformation,  we can reduce $\chi^{*}P$ to a model Klein-Gordon operator of the type introduced in Subsect.~\ref{sec1.1}. The reduction procedure is explained in details in \cite[Sect. 4]{GW1}. Theorem \ref{thm1} is then reduced to Thm.~\ref{thm1.2} below, whose proof will be explained in this section.
 \subsection{Model Klein-Gordon operators}\label{sec1.1}
Let us recall the model Klein-Gordon operators introduced in \cite[Sect. 2]{GW1}.
We work on $\rr^{1+d}$ with elements $x= (t, \rx)$ equipped 
with the Lorentzian metric
\[
g= - dt^{2}+ \alth_{ij}(t, \rx)d\rx^{i}d\rx^{j},
\]
where $t\mapsto \alth_{t}= \alth_{ij}(t, \rx)d\rx^{i}d\rx^{j}$ is a smooth family of Riemannian metrics on $\rr^{d}$.  Fixing a real potential $V\in \cinf(\rr^{1+d}; \rr)$ the Klein-Gordon operator $P = -\Box_{g}+V$ takes the form
\begin{equation}
\label{e1.0}
P= \p_{t}^{2}+ r(t, \rx)\p_{t}+ a(t, \rx, \p_{\rx}),
\end{equation}
where
\[
\begin{array}{l}
a(t)= a(t, \rx, \pe_{\rx})= - |\alth|^{-\12}\pe_{i}\alth^{ij}|\alth|^{\12}\pe_{j}+ \altV(t,\rx),\\[2mm]
 r(t)= r(t,\rx )= |\alth|^{-\12}\p_{t}(|\alth|^{\12})(t,\rx ).
\end{array}
\]
The operator $a(t)$ is formally selfadjoint for the  time-dependent scalar product
\[
(u|v)_{t}= \int_{\Sigma}\bar{u}v|h_{t}|^{\12}d\rx,
\]
and $P$ is formally selfadjoint for the scalar product
\[
(u|v)=  \int_{\rr\times \Sigma}\bar{u}v|h_{t}|^{\12}d\rx dt.
\]
Conditions $({\rm aM})$ on the original metric $\altg$ and potential $V$ imply similar asymptotic conditions on $a(t, \rx, \p_{\rx})$ and $r(t, \rx)$ when $t\to \pm \infty$. More precisely one has:
\[
({\rm Hstd}) \  \ \beal
a(t, \rx, \p_{\rx})= a_{\outin}(\rx, \p_{\rx})+ \Psi_{\std}^{2, - \delta}(\rr; \rr^{d})  \hbox{ on }\rr^{\pm}\times \rr^d,\\[2mm]
r(t)\in \Psi_{\std}^{0, -1-\delta}(\rr; \rr^{d}),\\[2mm]
a_{\outin}(\rx, \p_{\rx})\in \Psi_{\scc}^{2,0}(\rr^{d})\hbox{\ is  elliptic},\\[2mm]
a_{\outin}(\rx, \p_{\rx})= a_{\outin}(\rx, \p_{\rx})^{*}\geq C_{\infty}>0,
\eeal  
\]
where $\Psi_{\std}^{m, \delta}(\rr; \rr^{d})$ is the class of time-dependent pseudodifferential operators on $\rr^{d}$ associated to symbols $m(t, \rx, k)$ such that
\[
 \p_{t}^{\gamma}\p_{\rx}^{\alpha}\p_{\spexi}^{\beta}m(t, \rx, \spexi)\in O((\langle t\rangle+ \langle\rx\rangle)^{\delta- \gamma- |\alpha|}\langle \spexi\rangle^{m-|\beta|}),  \ \gamma\in \nn, \ \alpha, \beta\in \nn^{d}.
\]
Similarly $\Psi_{\scc}^{m, \delta}(\rr^{d})$ is the class of pseudodifferential operators on $\rr^{d}$ associated to symbols $m( \rx, k)$ such that
\[
\p_{\rx}^{\alpha}\p_{\spexi}^{\beta}m(\rx, \spexi)\in O( \langle\rx\rangle^{\delta-|\alpha|}\langle \spexi\rangle^{m-|\beta|}),  \ \alpha, \beta\in \nn^{d}.
\]
We refer the reader to \cite[Subsect.~2.3]{GW1} for more details.

\subsubsection{A further reduction}\label{sec1.1.1}
It is convenient to perform a further reduction to the case $r=0$. Namely setting $R= |h_{0}|^{\frac{1}{4}}|h_{t}|^{- \frac{1}{4}}$, we see that 
\[
L^{2}(\Sigma, |h_{0}|^{\12}d\rx)\ni\tilde{u}\mapsto R\tilde{u}\in L^{2}(\Sigma, |h_{t}|^{\12}d\rx)
\]
is unitary and that
\[
R^{-1}PR\eqdef\tilde{P}= \p_{t}^{2}+ \tilde{a}(t, \rx, \p_{\rx})
\]
where
\[
\tilde{a}(t)= rR^{-1}(\p_{t}R)+ R^{-1}(\p_{t}^{2}R)+ R^{-1}a(t)R
\]
is formally selfadjoint for $(\cdot| \cdot)_{0}$. Clearly $\tilde{a}(t, \rx, \p_{\rx})$ satisfies also $({\rm Hstd})$, with the same asymptotic $a_{\outin}(\rx, \p_{\rx})$. It is also immediate that the Hilbert spaces $\cY^{m}$, $\cX^{m}$, $\cX^{m}_{\F}$ introduced in Subsect.~\ref{sec1.3} are invariant under the map $u\mapsto Ru$ and hence  it suffices to prove Thm.~\ref{thm1.2} for $P$ replaced by $\tilde{P}$.

To simplify notation we will denote again $\tilde{P}$ by $P$. Summarizing we have $P= \p_{t}^{2}+ a(t, \rx, \p_{\rx})$, conditions $({\rm Hstd})$ are satisfied (with $r= 0$) and $a(t, \rx, \p_{\rx})$ is formally selfadjoint on
$L^{2}(\Sigma, |h_{0}|^{\12}d\rx)$ for all $t\in \rr$.

\subsection{Approximate diagonalization}\label{sec1.2}
Setting $\varrho u= \col{u}{\i^{-1}\p_{t}u}$, the equation $Pu=v$ is equivalent to
\[
(D_{t}- H(t))\varrho u= - \pi_{1}^{*}v, \ \ H(t)= \mat{0}{1}{a(t)}{0}
\]
where $\pi_{i}f= f_{i}$ for $f=\col{f_{0}}{f_{1}}$.

 One can  then construct an operator $T$ with $Tf(t)= T(t)f(t)$, and $t\mapsto T(t)$ is a smooth family of matrix-valued pseudodifferential operators on $\rr^{d}$ such that
\[
T^{-1}(D_{t}- H(t))T= D_{t}- H^{\rm ad}(t)\eqdef P^{\adg},
\]
where $H^{\rm ad}(t)$ is almost diagonal ie
\[
H^{\rm ad}(t)= H^{\rm d}(t)+ V^{\rm ad}_{-\infty}(t),
\]
\beq\label{e1.-1}
H^{\rm d}(t)= \mat{\epsilon^{+}(t)}{0}{0}{- \epsilon^{-}(t)}, 
\eeq
where $\epsilon^{\pm}(t)$ are time-dependent pseudodifferential operators on $\rr^{d}$, with principal symbols equal to $\pm (k\dual h^{-1}(t, \rx)k)^{\12}$, and $V^{\rm ad}_{-\infty}(t)$ is an off-diagonal matrix of time-dependent operators on $\rr^{d}$ such that 
 \[
 (\langle \rx\rangle + \langle t\rangle)^{m}V^{\rm ad}_{-\infty}(t)(\langle \rx\rangle + \langle t\rangle)^{-m+\delta}: H^{-p}(\rr^{d})\to H^{p}(\rr^{d})
 \]
  is uniformly bounded in $t$ for all $m, p\in \rr$.
Let us denote by $\cU^{\adg}(t, s)$ for $t,s\in \rr$ the Cauchy evolution generated by $H^{\adg}(t)$, i.e.~the solution of
\[
(D_{t}- H^{\adg}(t)) \cU^{\adg}(t,s)= 0, \ \ \cU^{\adg}(s, s)= \one.
\]
Then $\cU^{\adg}(t, s)$ is {\em symplectic}, which translates into the identity:
 \beq\label{e1.1}
H^{\rm ad}(t)^{*}q^{\rm ad}= q^{\rm ad}H^{\adg}(t), \ \ q^{\rm ad}\defeq\mat{1}{0}{0}{-1},
\eeq
where the adjoint is computed w.r.t.~the scalar product of $\cH^{0}= L^{2}(\rr^{d}, |h_{0}|^{\12}d\rx;  \cc^{2})$.

\subsection{Hilbert spaces}\label{sec1.3}
We set for $m\in \rr$
 \[
\cE^{\sobo}\defeq H^{\sobo+1}(\rr^{d})\oplus H^{\sobo}(\rr^{d}), \ \ \cH^{\sobo}\defeq H^{\sobo}(\rr^{d})\oplus H^{\sobo}(\rr^{d}), \ \sobo\in \rr,
\]
where $H^{m}(\rr^{d})$ are the usual Sobolev spaces. Fixing  $\gamma$ with $\12<\gamma<\12 + \delta$, we set:
\[
\cY^{\sobo}\defeq \langle t\rangle^{-\gamma}L^{2}(\rr; H^{\sobo}), \ \ \cY^{\adg, \sobo}\defeq   \langle t\rangle^{-\gamma}L^{2}(\rr; \cH^{\sobo}),
\]
and denote by $\cX^{m}$ the space of $u\in C^{1}(\rr; H^{m})\cap C^{0}(\rr; H^{m+1})$ such that $Pu\in \cY^{m}$, by $\cX^{\adg, \sobo}$ the space of $u^{\adg}\in C^{0}(\rr; \cH^{m})$ such that $P^{\adg}u^{\adg}\in \cY^{\adg}$. These spaces are equipped with the Hilbert space norms
 \beq\label{defdenorme}\begin{array}{l}
 \|u^{\adg}\|^{2}_{\sobo}\defeq \| \varrho^{\adg}_{0}u^{\adg}\|^{2}_{\cH^{\sobo}}+ \| P^{\adg}u^{\adg}\|^{2}_{\cY^{\adg, \sobo}},\\[2mm]
  \|u\|^{2}_{\sobo}\defeq \| \varrho_{0}u\|^{2}_{\cE^{\sobo}}+ \| Pu\|^{2}_{\cY^{\sobo}},
\end{array}
 \eeq
 where $\varrho_{t}u= \col{u(t)}{\i^{-1}\p_{t}u(t)}$, $\varrho^{\adg}_{t}u^{\adg}= u^{\adg}(t)$ are the Cauchy data maps at $t=0$. The well-posedness of the inhomogeneous Cauchy problems for $P$ and $P^{\adg}$, see \cite[Lemma 3.5]{GW1} implies that $\cX^{(\adg),m}$ are Hilbert spaces. 

\subsubsection{Feynman boundary conditions}\label{sec1.3.1}
Let us set
\[
\pi^{+}= \mat{1}{0}{0}{0}, \ \ \pi^{-}= \mat{0}{0}{0}{1}.
\]
It is straightforward to show that
\beq\label{e1.4}
\cX^{\adg, m}_{\rm F}\defeq\{u^{\adg}\in \cX^{\adg, m}: \lim_{t\to -\infty}\pi^{+}\varrho^{\adg}_{t}u^{\adg}= \lim_{t\to +\infty}\pi^{-}\varrho^{\adg}_{t}u^{\adg}=0\hbox{ in }\cH^{m}\}.
\eeq
is a closed subspace of $\cX^{\adg, m}$.  
In \cite[Sect. 3.7]{GW1} $\cX^{\adg, m}_{\rm F}$  is defined using scattering data maps, see \cite[Def.~3.7]{GW1}:  one sets 
\[
\epsilon_{\outin}= a_{\outin}(\rx, \p_{\rx})^{\12}, \ \ H_{\outin}^{\adg}= \mat{\epsilon_{\outin}}{0}{0}{-\epsilon_{\outin}}
\]
 and $\cU_{\outin}^{\adg}(t,s)= \e^{\i(t-s)H^{\adg}_{\outin}}$. Note that $H_{\outin}^{\adg}$ are actually diagonal \emph{exactly}, hence $\cU_{\outin}^{\adg}(t,s)$ commute with $\pi^+$ and $\pi^-$. The scattering data maps are 
 \[
\varrho^{\adg}_{\outin}u^{\adg}= \lim_{t\to \pm \infty}\cU_{\outin}^{\adg}(0, t)u^{\adg}(t)\hbox{ in }\cH^{m},
\]
in terms of which $\cX^{\adg, m}_{\rm F}$ can be defined as
\[
\cX^{\adg, m}_{\rm F}= \{u^{\adg}\in \cX^{\adg, m}: \pi^{+}\varrho_{\rm in}^{\adg}u^{\adg}= \pi^{-}\varrho_{\rm out}^{\adg}u^{\adg}=0\}.
\]
 Both definitions are the same using that $\cU^{\adg}(t, s)$ and $\cU^{\adg}_{\rm out/in}(t,s)$  are uniformly bounded in $B(\cH^{m})$, as shown in \cite[Prop.~5.6]{GW2}.

Similarly we define $\cX^{m}_{\F}$ as
\[
\cX^{m}_{\F}\defeq\{u\in \cX^{m}: \lim_{t\to-\infty}c^{-}_{\rm out}\varrho_{t}u= \lim_{t\to+\infty}c^{+}_{\rm in}\varrho_{t}u=0\hbox{ in }\cE^{m}\}
\]
where
$c^{\pm}_{\outin}= \12\mat{1}{\pm a_{\outin}^{\12}}{\pm a_{\outin}^{\12}}{1}$
 are the spectral projections on $\rr^{\pm}$ of $H_{\outin}= \mat{0}{1}{a_{\outin}}{0}$. Again, $\cX^{m}_{\F}$ is a closed subspace of $\cX^{m}$.

\subsection{Invertibility of $P$}
We now prove the following theorem:
\begin{theoreme}\label{thm1.2}
 $P: \cX^{m}_{\rm F}\to \cY^{m}$ is boundedly invertible with inverse
 \[
G_{\rm F}= - \pi_{0} T G^{\adg}_{\rm F}T^{-1}\pi_{1}^{*}.
\]
Moreover $G_{\rm F}$ is a {\em Feynman inverse} of $P$, i.e.,
\beq\label{e1.6}
 \WF(G_{{\rm F}})'= ({\rm diag}_{T^*M})\cup\textstyle\bigcup_{t\leq 0}(\Phi_t({\rm diag}_{T^*M})\cap \pi^{-1}\cN).
\eeq
\end{theoreme}

In view of the results in \cite{GW1}, the only part that deserves special attention is the proof that $P$ is injective, which is reduced to a similar statement about $P^{\adg}$.
The proof of Lemma \ref{lemma1.1} below is inspired by  the work of Vasy \cite[Prop.~7]{Va}, which in turn relies on  arguments of Isozaki \cite{I} from $N$-body scattering. In our framework it turns out to be very simple.
\begin{lemma}\label{lemma1.1}
 One has:
 \[
\Ker P^{\rm ad}|_{\cX^{\adg, m}_{\rm F}}= \{0\}\hbox{ for all }m\geq 1.
\]
\end{lemma}
\proof
 Let us set $\chi_{\epsilon}(t)= \int_{|t|}^{+\infty}\one_{[1, 2]}(\epsilon s)s^{-r}ds$ for some $0<r<1$. Note that $\supp \chi_{\epsilon}\subset \{|t|\leq 2 \epsilon^{-1}\}$. Let us still denote by $\chi_{\epsilon}$ the operator $\chi_{\epsilon}\otimes\one_{\cc^{2}}$.  Recalling that $q^{\adg}$ is defined in \eqref{e1.1},  we compute for $u^{\adg}\in \cX^{\adg, m}_{\rm F}$:
 \[
 \bea
&\int_{\rr}(P^{\adg}u^{\adg}(t)| q^{\adg}\chi_{\epsilon}(t)u^{\adg}(t))_{\cH^{0}}dt-(\chi_{\epsilon}(t)u^{\adg}(t)| q^{\adg} P^{\adg}u^{\adg}(t))_{\cH^{0}}dt\\[2mm]
&= \int_{\rr} (D_{t}u^{\adg}(t)| q^{\adg} \chi_{\epsilon}(t)u^{\adg}(t))_{\cH^{0}}- (u^{\adg}(t)| q^{\adg}\chi_{\epsilon}(t)D_{t}u^{\adg}(t))_{\cH^{0}}\\[2mm]
&\phantom{=}\, + \int_{\rr}  (u^{\adg}(t)| q^{\adg}[H^{\adg}(t), \chi_{\epsilon}(t)]u^{\adg}(t))_{\cH^{0}}dt,
\eea
\]
using that $H^{\adg *}(t)q^{\adg}= q^{\adg}H^{\adg}(t)$, $\chi_{\epsilon}(t)^{*}q^{\adg}= q^{\adg}\chi_{\epsilon}(t)$ and $u^{\adg}(t)\in \Dom H^{\adg}(t)$ since $m\geq 1$. 
We have $[H^{\adg}(t), \chi_{\epsilon}(t)]=0$,  and using  that $\chi_{\epsilon}$ is compactly supported in $t$ we can integrate by parts in $t$ in the second line and obtain
\begin{equation}
\label{e1.2}
\bea
&\int_{\rr}(P^{\adg}u^{\adg}(t)| q^{\adg}\chi_{\epsilon}(t)u^{\adg}(t))_{\cH^{0}}dt-(\chi_{\epsilon}(t)u^{\adg}(t)| q^{\adg} P^{\adg}u^{\adg}(t))_{\cH^{0}}dt\\[2mm]
&=-\i\int_{\rr}(u^{\adg}(t)| q^{\adg}\p_{t}\chi_{\epsilon}(t)u^{\adg}(t))_{\cH^{0}}dt.
\eea
\end{equation}
Note that we used here that the scalar product in $\cH^{0}$ does not depend on $t$, which is the reason for the reduction to $r=0$ in \ref{sec1.1.1}.

Since $P^{\adg}u^{\adg}=0$ this yields
\beq\label{e1.5}
\int_{\rr}(u^{\adg}(t)| q^{\adg}\p_{t}\chi_{\epsilon}(t)u^{\adg}(t))_{\cH^{0}}dt=0.
\eeq
We claim that:
\beq\label{e1.3}
\begin{array}{rl}
i)&\| \pi^{\pm}u^{\adg}(t)\|^{2}_{\cH^{0}}\in O(t^{1- \delta}) \hbox{ when }t\to \mp \infty,\\[2mm]
ii)&\| \pi^{\pm}u^{\adg}(t)\|^{2}_{\cH^{0}}= c^{\pm}+ O(t^{1-\delta})\hbox{ when }t\to \pm \infty.
\end{array}
\eeq
for $c^{\pm}= \lim_{t\to \pm\infty}\|\pi^{\pm}u^{\adg}(t)\|^{2}_{\cH^{0}}$. The proof of \eqref{e1.3}  is elementary:  we have $H^{\adg}(t)- H_{\outin}^{\adg}\in O(t^{-\delta})$ in $B(\cH^{0})$ when $t\to \pm \infty$, see e.g.~\cite[Subsect.~2.5]{GW1}, which  using $\delta>1$ and  the Cook argument (i.e., estimating first the derivative in time, and then integrating) yields:
\[
\begin{array}{l}
W_{\rm out/in}^{\dag}u^{\rm ad}= \lim_{t\to \pm\infty}\cU_{\rm out/in}^{\adg}(0, t)u^{\rm ad}(t)\hbox{ exists in }\cH^{0},\\[2mm]
\|W_{\rm out/in}^{\dag}u^{\rm ad}- \cU_{\rm out/in}^{\adg}(0, t)u^{\rm ad}(t)\|_{\cH^{0}}\in O(t^{1-\delta}).
\end{array}
\]
 Since $\cU_{\rm out/in}^{\adg}(0, t)$ is unitary on $\cH^{0}$ and $u^{\rm ad}(t)\in\cX^{m}_{\F}$, this implies \eqref{e1.3}. 

Next, we compute
\[
\bea
&\int_{\rr}(u^{\adg}(t)| q^{\adg}\p_{t}\chi_{\epsilon}(t)u^{\adg}(t))_{\cH^{0}}dt\\[2mm]
& = \int_{\rr}\p_{t}\chi_{\epsilon}(t)\|\pi^{+}u^{\rm ad}(t)\|_{\cH^{0}}^{2}dt- \int_{\rr}\p_{t}\chi_{\epsilon}(t)\|\pi^{-}u^{\rm ad}(t)\|_{\cH^{0}}^{2}dt.
\eea
\]
Since  $\p_{t}\chi_{\epsilon}(t)= -{\rm sgn}(t)\one_{[\epsilon^{-1}, 2\epsilon^{-1}]}(|t|)|t|^{-r}$ we have using \eqref{e1.3}:
\[
0\leq \int_{\rr^{\mp}}|\p_{t}\chi_{\epsilon}(t)|\|\pi^{\pm}u^{\rm ad}(t)\|_{\cH^{0}}^{2}dt\leq C\int\one_{[\epsilon^{-1}, 2\epsilon^{-1}]}(|t|) |t|^{-r- \delta+1}dt\in O(\epsilon^{r+ \delta-2}),
\]
and
\[
\bea
\int_{\rr^{\pm}}\p_{t}\chi_{\epsilon}(t)\|\pi^{\pm}u^{\rm ad}(t)\|_{\cH^{0}}^{2}dt&= \mp\int_{\rr^{\pm}} \one_{[\epsilon^{-1}, 2\epsilon^{-1}]}(|t|)c^{+}|t|^{-r} dt+ O(\epsilon^{r+ \delta-2})\\[2mm]
&=\mp Cc^{\pm}\epsilon^{r-1}+ O(\epsilon^{r+ \delta-2}).
\eea
\]
Using \eqref{e1.5} this yields $C\epsilon^{r-1}(c^{+}+ c^{-})\in O(\epsilon^{r+ \delta-2})$ hence $c^{+}= c^{-}=0$ since $\delta>1$. Therefore  by \eqref{e1.3} we have $\lim_{t\to \pm \infty}\| u^{\adg}(t)\|_{\cH^{0}}=0$. Since  the Cauchy evolution $\cU^{\adg}(t, s)$ is uniformly bounded in $B(\cH^{0})$ we have $u^{\adg}(0)=0$ hence $u=0$. \qed

\begin{proposition}\label{prop1.1}
  The operator $P^{\adg}: \cX_{\rm F}^{\adg, m}\to \cY^{\adg,m}$ is boundedly invertible for all $m\in \rr$.
  \end{proposition}
\proof
Recall that the  diagonal operator $H^{\d}(t)$  was introduced in \eqref{e1.-1}. In \cite[Thm.~3.16]{GW1} it is shown that $P^{\d}= D_{t}- H^{\d}(t):  \cX_{\rm F}^{\adg, m}\to \cY^{\adg,m}$ is boundedly invertible, with inverse $G^{\d}_{\rm F}$ given by:
 \[
\bea
G^{\d}_{\rm F} v^{\adg}(t)&\defeq  \i\int_{-\infty}^{t}\cU^{\d}(t, 0)\pi^{+}\cU^{\d}(0,s)v^{\adg}(s)ds\\ &\phantom{=\,}-\i\int_{t}^{+\infty}\cU^{\d}(t, 0)\pi^{-}\cU^{\d}(0,s)v^{\adg}(s)ds.
\eea
\]
Moreover by \cite[Lemma 3.7]{GW1} $V^{\adg}_{-\infty}$ is compact from $\cX^{\adg, m}$ to $\cY^{m}$  hence also from $\cX_{\rm F}^{\adg, m}$ to $\cY^{m}$ since $\cX_{\rm F}^{\adg, m}$ is closed in $\cX^{\adg, m}$. 

Therefore $P^{\adg}= P^{\d}+ V^{\adg}_{-\infty}: \cX_{\rm F}^{\adg, m}\to \cY^{m}$ is Fredholm of index $0$. If $u^{\adg}\in \Ker P^{\rm ad}|_{\cX^{\adg, m}_{\rm F}}$, we have $u^{\adg}= - G^{\d}_{\rm F}V^{\adg}_{-\infty}u^{\adg}$, from which we deduce that $u^{\adg}\in \cX^{\adg, m'}_{\rm F}$ for any $m'$, using that $V^{\adg}_{-\infty}$ is smoothing in $\rx$. By Lemma \ref{lemma1.1} we have  $\Ker P^{\rm ad}|_{\cX^{\adg, m}_{\rm F}}= \{0\}$ hence $P^{\adg}:\cX_{\rm F}^{\adg, m}\to \cY^{m}$  is boundedly invertible. \qed

\subsubsection{Proof of Thm.~\ref{thm1.2}}
 By \cite[(3.20)]{GW1} we know that  \[
 \pi_{0}T\in B(\cX^{\adg, \sobo+ \12}, \cX^{\sobo}), \ \ T^{-1}\pi_{1}^{*}\in B(\cY^{\sobo}, \cY^{\adg, \sobo+\12})
 \]
  hence $G_{\rm F}: \cY^{m}\to \cX^{m}$.  Since $(D_{t}- H(t))TG^{\adg}_{\rm F}T^{-1}=TG^{\adg}_{\rm F}T^{-1}(D_{t}- H(t))= \one$, we obtain that $P G_{\rm F}= G_{\rm F}P= \one$.
 We also have $\varrho\pi_{0}T G^{\adg}_{\rm F}T^{-1}\pi_{1}^{*}= T G^{\adg}_{\rm F}T^{-1}\pi_{1}^{*}v$.
From \cite[ equ. (3.25)]{GW1} we know that  $\varrho_{\overline{\rm F}}= \varrho^{\adg}_{\overline{\rm F}}T^{-1}\varrho$, hence $\varrho_{\overline{\rm F}}G_{\rm F}= 0$, i.e.~$G_{\rm F}: \cY^{m}\to \cX_{\rm F}^{m}$.  

To prove the second statement,  let  $\tilde{G}_{\rm F}= - \pi_{0} T G^{\d}_{\rm F}T^{-1}\pi_{1}^{*}$.  We have $G_{\rm F}^{\d}- G_{\rm F}^{\adg}\eqdef R= G^{\d}_{\rm F}V^{\adg}_{-\infty}G^{\adg}_{\F}$ by the resolvent identity.  It is shown in \cite[Lemma 3.7]{GW1} that $V^{\adg}_{-\infty}: \cX^{\adg, m}\to \cY^{m'}$ is bounded for all $m'>m$, hence $R: \cY^{\adg, m}\to \cX^{\adg, m'}$ for all $m'>m$ i.e.~it is smoothing in the $\rx$ variables. We use then that 
$D_{t}R= H^{\d}(t)R+ V^{\adg}_{-\infty}G^{\adg}_{\F}$, $RD_{t}= R H^{\adg}(t)+ G^{\d}_{\rm F}V^{\adg}_{-\infty}$ to gain regularity in the $t$ variable and obtain that $R: \cE'(\rr^{1+d}; \cc^{2})\to \cinf(\rr^{1+d}; \cc^{2})$.  Therefore, $G_{\rm F}- \tilde{G}_{\rm F}$ is a smoothing operator and it is shown in \cite[Thm.~3.18]{GW1} that $\WF(\tilde{G}_{{\rm F}})'$ equals the r.h.s.~of \eqref{e1.6}, which completes the proof. 
 \qed

\end{document}